\documentclass{svproc}
\usepackage{url}
\usepackage{graphicx}
\usepackage{smartdiagram}

\usepackage{amsmath}

\begin{document}
\mainmatter              
%
\title{On the Optimal Placement of Electric Vehicle Charging Stations}
%
%
\author{Johnny Tiu \and Shankar Ramharack \and Patrick Hosein}
\authorrunning{Johnny Tiu et al.} 
%
\tocauthor{Johnny Tiu, Shankar Ramharack, Patrick Hosein}
\institute{The University of the West Indies, St. Augustine, Trinidad and Tobago,
\email{johnnytiu.tt@gmail.com, shankar.ramharack@gmail.com, patrick.hosein@sta.uwi.com}}

\maketitle              
\begin{abstract} 
Increasing the adoption of Electric Vehicles (EV) is an integral part of many strategies to address climate change and air pollution. However, Electric Vehicle adoption rates are inhibited by several factors which reduce the confidence of potential buyers of electric vehicles, namely range anxiety and limited charging infrastructure. The latter concern can be addressed by carefully planning the placement of Electric Vehicle charging stations to sustainably meet long-term demand.
We outline a methodology to optimally allocate Electric Vehicle Charging stations using a novel approach to demand estimation. Using this approach the first two authors placed second in the Academic category of an International Shell competition which was created to obtain solutions to this problem. They were also provided with a grant to further study this approach and to develop a deployment strategy.

    
\end{abstract}

\section{Introduction}

Electric Vehicles (EVs) are distinguished by the use of electric motors as a primary driver for propulsion and inherently produce lower greenhouse gas (GHG) emissions compared to traditional Internal Combustion Engine (ICE) Vehicles. Thus, they play a major role in many countries' goals in achieving net-zero GHG emissions \cite{Paris}, with the United States of America setting a goal of 50\%  of all new vehicles sold in 2030 to be zero-emission \cite{whitehouse} and the United Kingdom designating that no diesel/gasoline-powered vehicles to be sold by 2030 \cite{uk-electric-vehicle-infr}. 

In addition, the purchase cost differential which is one of the major barriers to EV adoption is range anxiety  \cite{lee2018charging}, which is a fear that their EV will run out of charge before reaching their destination and is underpinned by the concern of available charging infrastructure. These each present their own issues but are ultimately intertwined and can be addressed through proper forecasting of electrical demand and optimal placement of charging stations.
Various modelling techniques have previously been used for short, medium and long-term electric demand forecasting such as structural time series and state space models, probabilistic models \cite{amini2016arima} , and Artificial Neural Networks (ANNs) \cite{jahangir2019charging} and Convolutional Neural Networks (CNNs) \cite{zhang2021deep}. We instead introduce a mixed-integer programming model for placement optimization together with a unique approach to demand estimation.

\section{Problem Description}

The Shell.ai Hackathon 2022 challenged teams to optimally place Electric Vehicle charging stations so that the configuration remains robust to demographic changes. They provided the required data, the constraints to be met and the objective to be achieved. Participants then had to develop a solution that would maximize the proposed objective function.

\subsection{Dataset Description}

\begin{figure}[ht]
\centering
    \includegraphics[width=0.6\textwidth]{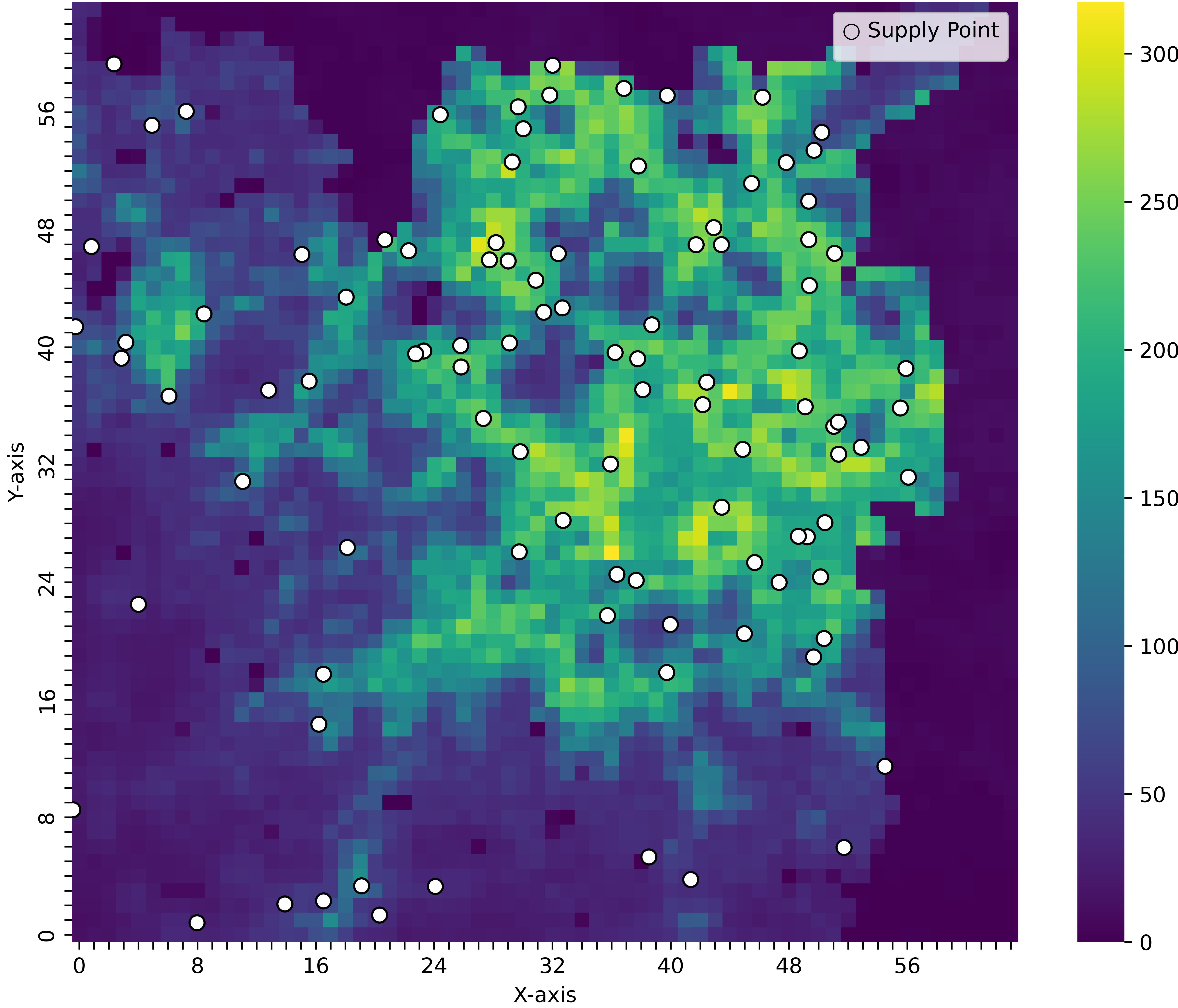}
    \caption{Charging Demand Heatmap and Supply point locations}
    \label{fig:enter-label}
\end{figure}

For an undisclosed location, a synthetic dimensionless dataset based on real world data was distributed to contestants. This consisted of:
\begin{description}
    \item[Demand History:] Yearly demand data between 2010 and 2018 for each demand point in a $64 \times 64$ grid.
    \item[Existing Infrastructure:] Distribution of Slow charging stations (SCS), Fast Charging Stations (FCS) and the maximum capacity for 100 supply points as of 2018. These are plotted in Figure \ref{fig:enter-label} showing the spatial distribution of supply points and charging demand intensity.
\end{description}

\subsection{Objective}

There are two objectives. Firstly, the EV charging station demand must be forecast over time and, using these predicted estimates, additional charging stations must be placed to satisfy this demand while minimizing the cost of the infrastructure.
The total cost comprised of three components:
\begin{description}
    \item[Customer Dissatisfaction:] This is a penalty based on how far customers have to travel to charge their vehicle and so represents charging anxiety. Let $d_{ij}$ denote the Euclidean distance between the $i^{th}$ demand point and the $j^{th}$ supply point and let $x_{ij}$ denote the amount of demand of the $i^{th}$ demand point that is satisfied by the $j^{th}$ supply point, then the customer dissatisfaction cost is given by
    \begin{equation}
        C_{cd} \equiv \alpha \sum_i \sum_j  x_{ij} d_{ij} 
    \end{equation}
    
    \item[Demand Mismatch:] This is the cost associated with incorrectly forecasting the EV charging demand. Let $D_i$ denote the true demand and $\hat{D}_i$ the predicted demand of demand point $i$ then 
    \begin{equation}
        C_{dm} \equiv \beta \sum_i |\hat{D}_i - D_i| \label{dm}
    \end{equation}
    \item[Infrastructure Cost:] This is the cost of operating, maintaining and developing the EV charging infrastructure. Let $N^{scs}_j$ denote the number of slow charging stations and $N^{fcs}_j$ the number of fast charging stations at the $j^{th}$ supply point then
    \begin{equation}
        C_{ic} \equiv \gamma \sum_j (N^{scs}_j + r N^{fcs}_j)  \label{if}
    \end{equation}
\end{description}
The scaling parameters were given as $\alpha=1$, $\beta=25$, $\gamma=600$ and $r=1.5$ which is the ratio of costs associated with developing Fast and Slow charging stations.

\subsection{Problem Constraints} \label{Constraints}
The constraints we provided as follows:
\begin{enumerate}
    \item All elements of the demand supply matrix $\mathbf{X}$ must be non-negative.
    \item The number of slow $N^{scs}_j$ and fast $N^{fcs}_j$ chargers must be non-negative.
    \item The total number of chargers placed at a supply point $j$ must not exceed the number of parking spots which we denote by $p_j$. 
    \item Charging stations can only be added (not removed) over time at a supply point.
    \item The total demand satisfied at a supply point cannot exceed the available supply.
    \item The predicted demand at each demand point must be exactly satisfied.
\end{enumerate}

\section{Mathematical Models}
The two aspects of the problem, demand estimation and placement optimization, are independent and so can be solved independently. We provide the models that were used for each and outline the solution approach. 


\subsection{EV Demand Prediction Model} \label{Prediction Model}

For each demand point we are given actual demand over a number of years and we need to predict the demand for the next two years.
For each demand point we can use the data from 2010 to 2018 to predict the value for 2019. The issue was that the demand for a particular point may vary widely from year to year which resulted in a poor prediction. We therefore decided to used points surrounding the demand point to get a better estimate of its demand. For a given year (with known demand) we instead estimate the demand for point $k$ as follows:
\begin{equation}
    \tilde{D}_k = \frac{\sum_j \frac{D_j}{(1 + {d_{jk})^{\kappa}}}}{\sum_j \frac{1}{1 + d_{jk})^\kappa}}
\end{equation}
Note that when $\kappa=0$ we use the average demand over all demand points while if $\kappa=\infty$ then $\tilde{D}_k = D_k$. This approach was adapted from an approach used in \cite{hosein2023data}.

These modified demands $\tilde{D}_j$ are computed for each year (2010-2017) for a given value of $\kappa$ and then a third order polynomial regression algorithm is used to predict demand for 2018. We then compared the true and predicted demands for 2018 to obtain the Mean Square Error (MSE) for this value of $\kappa$. This is repeated for different values of $\kappa$. We found that $\kappa=4.7$ provided the lowest MSE. Finally we used this value of $\kappa$ for years 2010-2018 to predict the demands for 2019 and then for 2020. These demands were then used in the charging station placement problem which we next describe.

\subsection{Charging Station Placement Optimisation Model}

\begin{figure}[ht]
\centering
\input{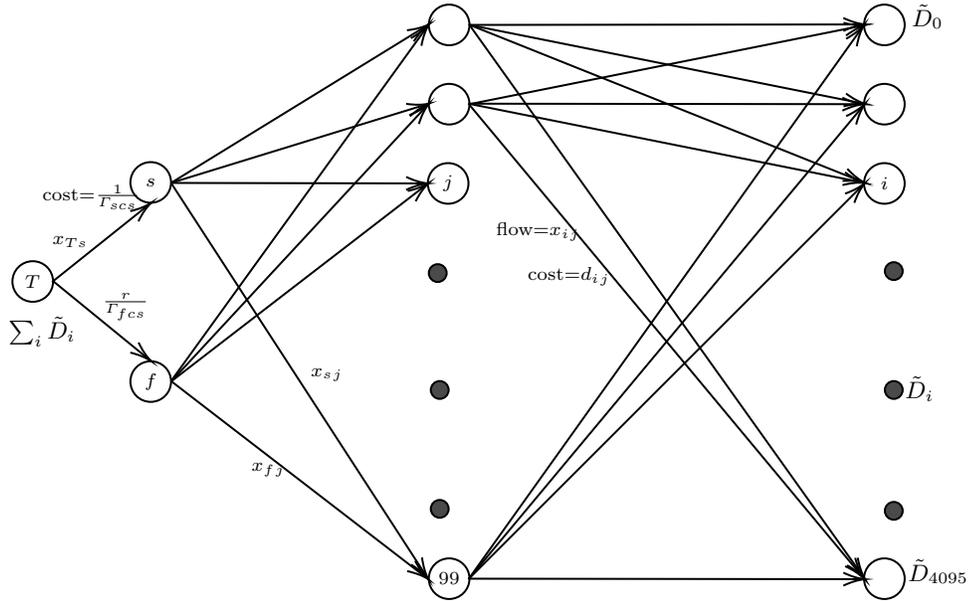}
\caption{Network Flow Model}
\label{net}
\end{figure}

The optimal placement of chargers can be formulated as a Mixed-Integer Programming Problem.
Note that, since the true demands for 2019 and 2020 are not known then we cannot include this cost in our model (although the cost was used by the organizers when evaluating the solution. The optimization problem and constraints can therefore be stated as:

\begin{equation}
\min_{\mathbf{X}, \vec{N^{scs}}, \vec{N^{fcs}}} Z = \alpha \sum_{i,j} (x_{ij} d_{ij}) + \gamma \sum_j (N^{scs}_j + r N^{fcs}_j)
\end{equation}
Subject To:
\begin{equation}
x_{ij} \ge 0 \quad \forall \; i, j
\end{equation}
\begin{equation}
N^{scs}_j \in \{0, 1, 2, \dots\}, \quad N^{fcs}_j \in \{0, 1, 2, \dots\} \quad \forall \; j
\end{equation}
\begin{equation}
N^{scs}_j + N^{fcs}_j \leq p_j \;\;\; \forall j
\end{equation}
\begin{equation}
\sum_i x_{ij} \leq \Gamma_{scs} N^{scs}_j + \Gamma_{fcs} N^{fscs}_j \quad \forall \; j
\end{equation}
\begin{equation}
N^{scs}_j \geq N^{scs}_j(2018) \;\;\; N^{fcs}_j \geq N^{fcs}_j(2018)
\end{equation}
\begin{equation}
\sum_j x_{ij}  = \hat{D}_i \quad \forall \; i
\end{equation}
where $r = 1.5, \alpha = 1$ and $\gamma=600$. The capacity of a slow charging station is $\Gamma_{scs}=200$ and that of a fast charging station is $\Gamma_{fcs}=400$. Figure \ref{net} contains a network flow formulation of the problem.

This minimisation problem was solved using the Mixed-Integer Programming Solver in Julia using the Branch and Cut option and using the projected demands obtained from Section \ref{Prediction Model} and data on the existing infrastructure. A lower bound to this problem can be obtained by relaxing the integer constraints and solving the resulting Linear Programming problem exactly. The Mixed-Integer Programming solution (which may not be optimal) was within 0.01\% of the lower bound which means it was near-optimal.

\section{Results and Leader board ranking} 

The information required to evaluate the demand prediction was only available to the organizers and hence they alone could compute the demand mismatch cost. However we have shown that, given demands, the placement optimization solution is already very close to optimal (within 0.01\% within the optimal objective function value).

In terms of the demand estimation aspect. The prediction produced a result that performed second in the academic category (meaning that one other University had a better solution). Our approach outperformed methods employing polynomial regressions combined with MILP, ANN, CNN and ARIMA models.

\section{Conclusion and Future Work}
We outlined an approach for determining the optimal placement of Electrical Vehicle charging stations. This consisted of two subproblems. We first estimated demand for charging stations based on historical data. We introduced a novel approach for estimating demand that allows one to improve the robustness of the estimate by using nearby samples. We formulated the charging station placement optimization problem as a Mixed Integer Programming problem and demonstrated that the solution obtained was within 0.01\% of the optimal solution. We are in the process of improving these algorithms using a grant received from Shell. In the future we will report on any new results.

\bibliographystyle{spmpsci_unsrt}
\bibliography{evdemand}
\end{document}